\documentclass[smallextended,final]{svjour3}

\usepackage{graphics,graphicx}
\usepackage{upquote}
\usepackage{verbatim}
\newlength\myverbindent 
\def\C{{\bf C}}
\def\Re{\kern .7pt\hbox{Re}\kern .7pt}
\def\Im{\kern .7pt\hbox{Im}\kern .7pt}
\def\Omegab{\overline{\Omega}}
\def\Db{\overline{D}}
\def\pput(#1,#2)#3{\noindent\smash{\raise#2pt\hbox to 0pt
   {\kern #1pt #3\hss}}\ignorespaces}

\title{Numerical conformal mapping with rational functions}

\author{Lloyd N.~Trefethen}
\institute{Lloyd N. Trefethen \at Mathematical
Institute, University of Oxford, Oxford, OX2 6GG, UK,
\email{trefethen@maths.ox.ac.uk}}

\begin{document}

\maketitle

\begin{abstract}
New algorithms are presented for numerical conformal mapping based
on rational approximations and the solution of Dirichlet problems by
least-squares fitting on the boundary.  The methods are targeted at
regions with corners, where the Dirichlet problem is solved by the
``lightning Laplace solver'' with poles exponentially clustered
near each singularity.  For polygons and circular polygons,
further simplifications are possible.  
\end{abstract}

\section{\label{sec-intro}Introduction}
The association of conformal mapping with Laplace problems goes back
to the beginning of the subject~\cite{henrici,schiffer}.  When this
connection has been exploited for the numerical computation of
conformal maps, it has usually been by way of integral equations
such as those of Berrut, Gerschgorin, Lichtenstein, Symm, and
Warschawski~\cite{gaier,henrici,tref,wegmann}.  However, almost all
domains of interest in practical conformal mapping contain corners,
and for such domains, it has recently been discovered that rational
functions with exponentially clustered poles offer a powerful
alternative approach to solving Laplace problems~\cite{lightning}.
The aim of this paper is to bring these ideas together to present
simple new algorithms for numerical conformal mapping.  A second,
different aspect of rational functions introduced in~\cite{gt} also
comes into play: the compression of a computed conformal map to a
more compact form by means of the AAA algorithm~\cite{AAA,lawson}.

As shown in~\cite{lightning}, rational functions can achieve
root-exponential convergence for almost arbitrary corner
singularities, i.e., $O(C^{\kern .5pt -\sqrt n}\kern .4pt )$ as a
function of the rational degree~$n$ for a constant $C>1$, and no
analysis of the singularities is needed.  This means that our method
does not care if the boundary components are straight or curved, so
long as they are smooth.  Unbounded domains are also not a problem.

For simply-connected domains with smooth boundaries, the idea of
conformal mapping via Laplace problems can be realized by means of
polynomials instead of rational functions, and we begin in Section~2
by presenting an algorithm of this kind.  Though this method is
sometimes effective, it often fails for reasons related to the
``crowding phenomenon'' of numerical conformal mapping.  Methods
based on integral equations are superior, and such a method is the
default used in the new Chebfun code {\tt conformal}.  The picture
changes sharply when we turn to rational functions in Section~3.  We
show that a conformal map of a simply-connected domain with corners
can be computed in a few lines of code by calling the \hbox{MATLAB}
``lightning Laplace solver'' {\tt laplace}~\cite{lightcode}.

Both the polynomial and rational function methods are readily
generalized to doubly connected domains, and these generalizations
are presented for smooth boundaries in Section~\ref{smooth2} and for
boundaries with corners in Section~\ref{corners2}.  Domains of higher
connectivity are considered in Section~6, and numerical details of
our methods are discussed in Section~7.  Section 8 discusses the
further simplifications available for polygons and circular polygons
if one uses a modified method based on explicit information about
corner singularities.  Conclusions are presented in Section~9.

\section{\label{smooth1}Smooth simply connected domains}
Let $\Omega\in \C$ be a simply connected domain bounded by a
Jordan curve $P$, and let $D$ and $S$ be the open unit disk and the unit circle.
Without loss of generality we assume $0\in \Omega$ and
seek the unique conformal map $f: \Omega \mapsto D$ such that
$f(0) = 0$ and $f'(0)>0$.  By the Osgood--Carath\'eodory
theorem~\cite{henrici}, $f$ extends to a homeomorphism
of $\Omegab$ onto $\Db$.  It follows that $g(z) = \log(f(z)/z)$ is a nonzero
analytic function on $\Omega$ that is continuous on $\Omegab$ and
has real part $-\log |z|$ for $z\in P$ and imaginary part
$0$ at $z=0$.  If we write $g$ as
$g(z) = u(z) + iv(z)$, where $u$ and $v$ are real harmonic functions
in $\Omega$, continuous on $\Omegab$, then $u$ is the unique solution of
the Dirichlet problem
\begin{equation}
\Delta u = 0\, ; \quad u(z) = -\log|z|, ~ z\in P,
\label{dirichlet}
\end{equation}
and $v$ is its unique harmonic conjugate in $\Omega$ with $v(0) = 0$.
Combining these elements, we see that $f$ takes the form
\begin{equation}
f(z) = z e^{u(z) + i v(z)}.
\label{frep}
\end{equation}
Note that $u(z) + \log |z|$ is the Green's function of $\Omega$
with respect to the point~$z=0$, so $f$ is nothing more than the
exponential of the analytic extension of the Green's function.
See~\cite[Thm.~16.5a]{henrici} and~\cite[p.~253]{schiffer}.

To compute the conformal map $f$, therefore, it is enough to solve
the Dirichlet problem (\ref{dirichlet}) for $u$ and then find its
harmonic conjugate $v$.  In this paper we approximate solutions to
Dirichlet problems by the real parts of analytic functions defined
by series involving polynomials, rational functions, or other
elementary terms.  The harmonic conjugates come for free from the
same expansion coefficients.

In this section we consider the simplest case, in which $P$ is
well-behaved enough that it is effective to work with polynomials.
We first outline the algorithm in the monomial basis and then
mention the variation involving Arnoldi factorization
that our codes actually employ.
Following~\cite{series}, let $u$ be approximated as
\begin{equation}
u(z) \approx \Re\kern -1pt  \sum_{k=0}^n c_k z^k
\label{approx0}
\end{equation}
for some degree $n$
and complex coefficients $c_k$, which will depend on $n$.
Writing $c_k = a_k + ib_k$, we have
\begin{equation}
u(z) \approx \sum_{k=0}^n \,\bigl( a_k \Re(z^k) - b_k \Im(z^k)\bigr).
\label{approx}
\end{equation}
With $b_0=0$, this gives us $N = 2n+1$ real unknowns
$a_0,\dots, a_n$ and $b_1,\dots, b_n$.  These are determined
numerically by sampling the function $-\log|z|$ in $M\gg N$ points
on $P$ and solving an $M\times N$ matrix least-squares problem
of the form $Ax\approx b$\kern .7pt; we take $M\approx 4N$.
The harmonic conjugate is then approximated by the series
\begin{equation}
v(z) \approx \sum_{k=0}^n \,\bigl(a_k \Im(z^k) + b_k \Re(z^k)\bigr),
\label{vapprox}
\end{equation}
and $f$ comes from (\ref{frep}).  Note that if $P$ is analytic, then
$f$ can be analytically continued to a neighborhood of $\Omegab$,
and geometric convergence of (\ref{approx}) and (\ref{vapprox})
is guaranteed (in exact arithmetic) as $n\to\infty$~\cite{walsh},
assuming the boundary sampling is fine enough that the discrete
least-squares problem on $P$ approximates the continuous one.
(There are also special cases in which $P$ is not smooth and
yet $f$ is still analytic everywhere on $P$, such as a rectangle,
an equilateral triangle, or the ``circular L-shape'' consisting of
a square from which a quarter-circular bite has been removed.)

By the method just described, if $\Omega$ is well enough behaved,
one can compute a polynomial-based approximation of the conformal
map $f$.  The degree $n$ will usually have to be high.  However,
as described in~\cite{gt}, the computed map can be compressed to a
much more compact form as a rational function by means of the AAA
approximation algorithm~\cite{AAA}.  Moreover, the AAA method can
be invoked a second time to compute a rational approximation of
the inverse map $f^{-1}$.  As emphasized in~\cite{gt}, when AAA
approximation is brought into play, the traditional fundamental
difference between conformal maps in the forward and inverse
directions ceases to be very important.  Once one is available,
the other is readily computed too.

The presentation of the last page has been framed with respect
to monomial bases $1, z, z^2,\dots.$  Unless $P$ is close
to a circle about the origin, however, the monomials will be
very ill-conditioned, and it is better to use an approximation
to a family of polynomials orthogonal over $P$.  This we do by
carrying out an Arnoldi iteration as described in~\cite{VmeetsA}.
Mathematically the two approaches are equivalent.

We summarize the algorithm, which we call {\em Algorithm P1\/}
(P for polynomial, 1 for simply connected), as follows.  This is a
high-level description, omitting many details, and the parameters
indicated are engineering choices, not claimed to be optimal.
For reasons discussed in~\cite{gt}, high accuracy tends to be
problematic in numerical conformal mapping, and we work by default
to a tolerance of $10^{-5}$.

\medskip
{\em\openup 2pt

Algorithm P1.

$\bullet$ For $n = $ rounded values
of\/ $2^4, 2^{4.5}, 2^5,\dots$ until convergence or failure,

\indent\kern .19in - Set $M=8\kern .3pt n$ and sample
${-}\kern 1pt\log|z|$ in $M$ points of\/ $P$;

\indent\kern .19in - Solve an $M\times N$ least-squares
problem for the approximation $(\ref{approx})$;

\indent\kern .19in - For stability, do this via
Arnoldi factorization as described in\/~{\rm\cite{VmeetsA}}.

$\bullet$ Form $v$ from $(\ref{vapprox})$ and $f$ from
$(\ref{frep})$\kern .7pt;

$\bullet$ Approximate $f$ by a AAA rational approximation\kern .7pt;

$\bullet$ Approximate $f^{-1}$ by another AAA rational approximation.

\par}

\medskip

\noindent This algorithm is available in Chebfun~\cite{chebfun}
in the command \kern .4pt{\tt conformal}, introduced in October,
2019, if it is called with the flag \verb|'poly'|.  It is not the
default algorithm used by \kern .4pt{\tt conformal}, as we shall
explain in a moment.

\begin{figure}
\begin{center}
\vskip .15in
\includegraphics[scale=.67]{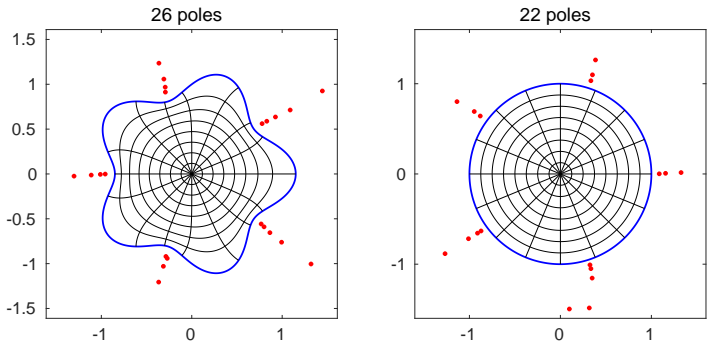}
\end{center}
\caption{\label{snowflake}Conformal map of a
smooth domain computed to about\/ $5$-digit accuracy
in\/ $0.2$ seconds on a desktop machine by
Algorithm \hbox{\rm P1} via the 
Chebfun \kern .4pt{\tt conformal(...,'poly')} command.
The red dots mark poles of the rational
representations of\/ $f$ and $f^{-1}$: just\/ $26$ poles for $f$
and\/ $22$ for $f^{-1}$. 
The closest poles lie at distance $0.099$ from $P$ on
the left and at distance $0.059$ from $S$ on the right.}
\end{figure}

Figure~\ref{snowflake} illustrates Algorithm P1 applied to a smooth
domain with five-fold symmetry, showing the output from the Chebfun
commands

{
\begin{verbatim}

     C = chebfun('exp(1i*pi*t)*(1+.15*cos(5*pi*t))', 'trig');
     [f, finv] = conformal(C, 'poly', 'plots');
\end{verbatim}
\par}

\noindent
The objects \kern.6pt{\tt f} and \kern.4pt{\tt finv} are function
handles corresponding to rational approximations of degrees 26 and
22, respectively, of the functions $f$ and $f^{-1}$.  These degrees
are remarkably low.  The computation takes about $0.3$ seconds on a
desktop machine, plus another second for plotting, and the maximal
error on the boundary is about $10^{-5}$.

\begin{figure}
\begin{center}
\includegraphics[scale=.67]{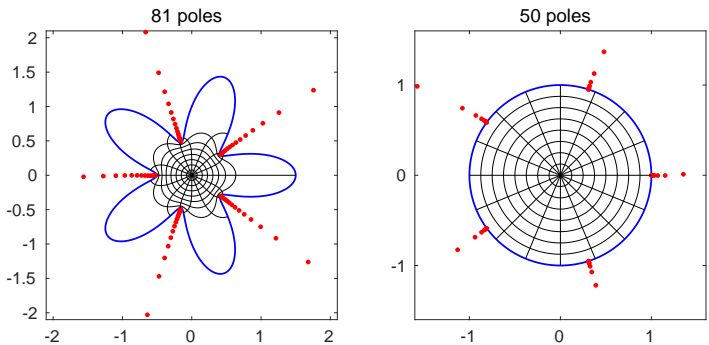}
\end{center}
\caption{\label{snowflakedeep}If the lobes are made deeper
by changing the coefficient of the cosine term from $0.15$ to $0.5$, Algorithm~P1
fails because it would require a polynomial of degree in the
tens of thousands.  Instead this figure was computed
via the Kerzman--Stein integral equation, the default algorithm
of Chebfun~\kern .4pt{\tt conformal}.
The closest poles now lie at distances $0.013$ from $P$ on
the left and $0.0011$ from $S$ on the right.
The AAA rational representations are still remarkably economical.}
\end{figure}

All this sounds good, but in fact, Algorithm P1 is limited in its
capabilities, and will fail for many smooth domains.  For example,
it fails if the indentations in the domain of Figure~\ref{snowflake}
are made twice as deep by changing the coefficient $0.15$ to $0.30$.
The problem in such a case is that the polynomial degree would have
to be in the thousands, a consequence of the large distortions
in conformal mapping that go by the name of the {\em crowding
phenomenon}~\cite{wegmann}.  If $0.30$ is increased further to
$0.50$, we get the domain shown in Figure~\ref{snowflakedeep}, and
this was certainly not mapped with Algorithm~P1.  The default method
of {\tt conformal} was applied instead, involving a discretization
of the Kerzman--Stein integral equation~\cite{ks,kt} followed by
AAA rational compression in both directions.  For this domain,
Algorithm~P1 would require a polynomial of degree in the tens
of thousands.

There is a literature on the crowding phenomenon, and in particular,
see Theorems 2--5 of~\cite{gt}.  These establish that if $\Omega$
contains an outward-pointing finger of length-to-width ratio $L$
as defined in a certain precise sense, then $f^{-1}$ must be at
least as large as order $e^{\pi L}$ at some points of $S$, with
singularities (or at least points of non-univalency) at distances
no greater than order $e^{-\pi L}$ from $S$, and that any useful
polynomial approximation $p\approx f^{-1}$ must have degree at least
on the order of $e^{\pi L}$.  (In Figure~\ref{snowflakedeep} we have
$L\approx 1.75$ and $e^{\pi L}\approx 244$.)  Analogous if somewhat
weaker effects occur for $f$ and its polynomial approximations
if there are inward-pointing fingers.  The geometric convergence
mentioned earlier still applies, but is useless in practice
because the constant $C$ exceeds $1$ by only a small amount.
These effects are reflected in the positions of the red dots
marking poles in Figures~\ref{snowflake} and~\ref{snowflakedeep},
highlighting the profound difference between polynomial and rational
approximations.  In the right image of each pair, the poles cluster
near points of $S$ that map to outward-pointing lobes of $\Omega$.
Similarly, the poles in the left image cluster near points of $P$
at inward-pointing fingers.

We emphasize that the effects just discussed will cause trouble
for any method based on polynomial approximations, no matter how
stably implemented, for they have nothing to do with the choice of
polynomial basis.  Our Arnoldi factorization makes the
implementation quite stable, but still it will fail for
domains whose maps involve large distortions.

Variants of Algorithm P1 can be used to solve related problems
involving exterior domains.  To map the interior of $P$ to
the exterior of $S$, we change the factor $z$ to $z^{-1}$
in (\ref{frep}); the Dirichlet problem (\ref{dirichlet})
remains unchanged.  To map the exterior of $P$ to the
exterior of $S$, we replace positive by negative powers in
(\ref{approx0})--(\ref{vapprox}).  To map the exterior of $P$
to the interior of $S$, we make both of these changes.

We now turn now to the much more powerful Algorithm~R1 based on
rational functions.

\section{\label{corners1}Simply connected domains with corners}

Let the boundary $P$ of~$\Omega$ be a polygon or circular polygon or
other Jordan curve consisting of a finite set of analytic curves
meeting at corners $z_1,\dots,z_K$.  (A circular polygon is a
finite chain of circular arcs.)  As before, we assume $0\in \Omega$
and seek the conformal map $f: \Omega \mapsto D$ with $f(0) = 0$
and $f'(0)>0$.

Equations (\ref{dirichlet}) and (\ref{frep}) still hold:
we need to solve a Dirichlet problem, but now, there will be
singularities at the corners.
Instead of (\ref{approx0}), we consider
\begin{equation}
u(z) \approx \Re \Biggl[\;
\sum_{j=1}^{n_1}{a_j\over z-z_j}
\;+\; \sum_{j=0}^{n_2}b_j z^j\Biggr]
\label{form}
\end{equation}
with complex coefficients $a_j$ and $b_j$; again, the polynomial
part of the sum will be implemented in practice by an Arnoldi
factorization~\cite{VmeetsA}.  Since the
publication of a paper by Newman in 1964~\cite{newman}, it
has been recognized that rational functions of the form of the
first sum in brackets can approximate certain singularities with
root-exponential convergence, if the poles~$z_j$ are clustered
exponentially near the singular points with a spacing that scales
with $n_1^{-1/2}$.  In~\cite{lightning} it was proved that this
effect applies generally to solutions of Lapalce
Dirichlet problems defined by piecewise analytic functions meeting
at corners.  In particular, Theorem~2.3
of~\cite{lightning} implies that approximations of the form
(\ref{form}) can converge root-exponentially for the conformal
map~$f$.  (The proof given in~\cite{lightning} assumes $P$ is
convex, but it is believed that the result holds generally.)
Since the appearance of~\cite{lightning}, a MATLAB ``lightning
Laplace solver'' code {\tt laplace} has been developed to solve
Dirichlet problems by this method~\cite{lightcode}.

We summarize our algorithm, {\em Algorithm R1}.

\medskip
{\em\openup 2pt

Algorithm R1.

$\bullet$ Solve problem~$(\ref{dirichlet})$ by
the lightning method and form $f$ from {\rm (\ref{frep})};

$\bullet$ Approximate $f$ by a AAA rational approximation\kern .7pt;

$\bullet$ Approximate $f^{-1}$ by another AAA rational approximation.

\par}

\medskip

\noindent Using \kern .4pt{\tt laplace}\kern .7pt, all
this can be done in four lines of MATLAB, assuming a specification
of $P$ and a tolerance {\tt tol} have been
given.  If $P$ is just a polygon,
then {\tt P} can be just a vector of vertices.
By default we take $\hbox{\tt tol} = 10^{-6}$:

{
\begin{verbatim}

     [~, ~, w, Z] = laplace(P, @(z) -log(abs(z)), 'tol', tol);
     W = Z.*exp(w(Z));
     f = aaa(W, Z, 'tol', tol);
     finv = aaa(Z, W, 'tol', tol);
\end{verbatim}
\par}

\noindent This is the central message of this paper: the conformal
map of a domain with corners can be computed in this simple fashion.
Section~7 will discuss certain numerical details.

\begin{figure}
\begin{center}
\vskip .15in
\includegraphics[scale=.67]{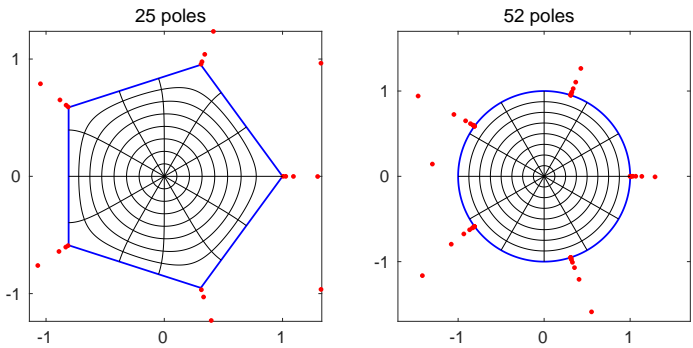}
\end{center}
\caption{\label{pentagon}Conformal map of a
pentagon computed in about\/ $0.3$ secs.\ by
Algorithm \hbox{\rm R1}.}
\end{figure}

Figure~\ref{pentagon} illustrates Algorithm R1 applied to
a regular pentagon, a polygonal analogue of the domain of
Figure~\ref{snowflake}.  Five-digit accuracy is achieved in
about $0.3$ secs.\ of desktop time, with poles in both directions
clustering near the five singular points.  There are four poles
near the vertex at $z=1$, for example, whose distances from the
vertex are approximately

{
\begin{verbatim}

     0.0097,  0.044,  0.14,  0.49.
\end{verbatim}
\par}
\noindent Though these poles were chosen adaptively by the AAA
algorithm, based on no mathematical analysis, they show exponential
clustering with a ratio of successive distances on the order of $3$.

\begin{figure}
\vskip .2in
\begin{center}
\includegraphics[scale=.67]{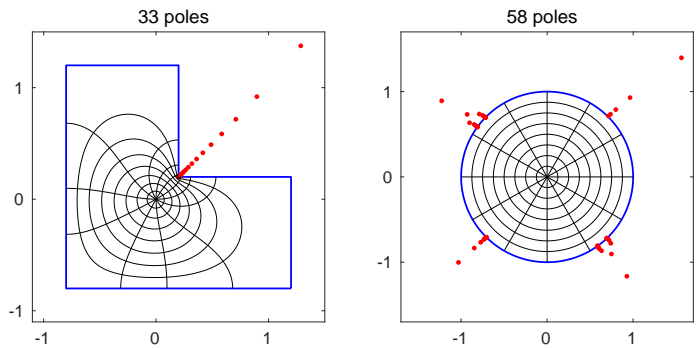}
\end{center}
\caption{\label{L}Mapping of an L-shaped domain by Algorithm R1.
The map $f$ is nonsingular
at the salient corners, so on the left, poles cluster
only at the reentrant corner, whereas on the right
they cluster at six points.}
\end{figure}

\begin{figure}
\vskip .2in
\begin{center}
\includegraphics[scale=.35]{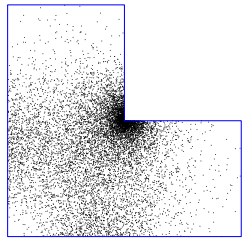}
\end{center}
\caption{\label{Ldots}The rational approximations produced
by Algorithm R1 can be evaluated in just microseconds.
This figure shows images under $f^{-1}$
of\/ {\rm 10,000} uniformly distributed random points
in the unit disk, computed in about $0.005$ secs.\ total.  The
nonuniform shading illustrates the exponential distortion effects in
conformal maps, discussed in Section\/~$2$.}
\end{figure}

Figure~\ref{L} shows results of an analogous computation, taking
about half a second, for an L-shaped domain.  On the left, poles
cluster near the reentrant corner at distances approximately

{
\begin{verbatim}

    0.000053, 0.00030, 0.00069, 0.0071, 0.0013, 0.0022, 0.0035,
    0.0054, 0.081, 0.012, 0.017, 0.024, 0.034, 0.047, 0.066,
    0.091, 0.12, 0.17, 0.23, 0.31, 0.42, 0.56, 0.76, 1.04, 1.66.
\end{verbatim}
}
\noindent These show exponential spacing with a factor of about 1.4.
The reduction of the factor from 3 to 1.4 as the number of
poles near a corner increases from 4 to 25 is consistent with
root-exponential clustering theorems of~\cite{lightning}.

We emphasize the extraordinary economy of the representations of
conformal maps computed by Algorithm R1, since they are nothing
more than rational functions of a degree typically less than~100.\
\ Figure~\ref{Ldots} shows images of $10^4$ points in the unit disk,
computed in $0.005$ secs.\ all together.

\begin{figure}
\vskip .2in
\begin{center}
\includegraphics[scale=.67]{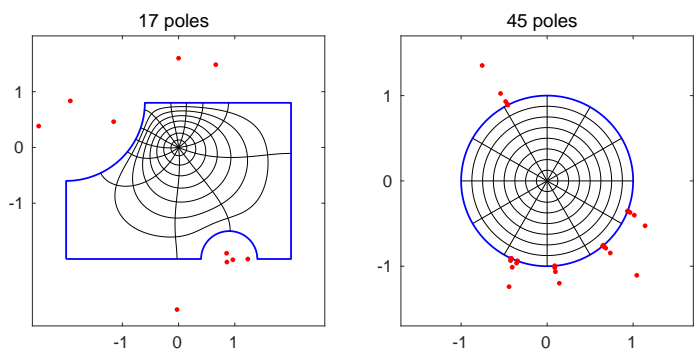}
\end{center}
\caption{\label{howell}Conformal map of a circular polygon
computed by Algorithm~{\rm R1}.  This shape is the first
example presented by Howell in~{\rm\cite{howell}}.}
\end{figure}

\begin{figure}
\vskip .2in
\begin{center}
\includegraphics[scale=.67]{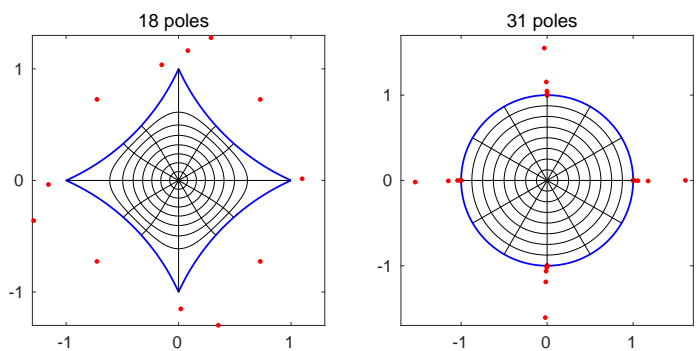}
\end{center}
\caption{\label{scalloped}Another map of a circular polygon computed
by Algorithm R1.  Each boundary component is a circular arc of radius $2$.}
\end{figure}

These examples involve polygons, but Algorithm R1 is equally applicable to
other domains bounded by analytic arcs
meeting at corners.  For example, a {\em circular polygon} is a generalization
of a polygon in which each side may be a straight segment or a circular arc.  Whereas
polygons are conventionally mapped by the Schwarz--Christoffel integral,
circular polygons can be mapped by the so-called Schwarzian differential equation~\cite{SCbook}.
This idea goes goes back to Schwarz himself 150 years ago, but has rarely been
attempted numerically; the two realizations we know of are~\cite{bg} and~\cite{howell}.
Figures~\ref{howell} and~\ref{scalloped} show maps of circular
polygons computed by the completely different approach of
Algorithm R1, each in less than a second on a desktop.  The {\tt laplace} code
used for these computations includes a convenient syntax for specifying
circular boundary arcs via their initial point and radius of curvature. 
Thus the mapping of Figure~\ref{scalloped} was computed by the
four-line command sequence listed above preceded by the specification
\begin{verbatim}

     P = { [1 -2], [1i -2], [-1 -2], [-1i -2] };
\end{verbatim}

\begin{figure}
\vskip .2in
\begin{center}
\includegraphics[scale=.67]{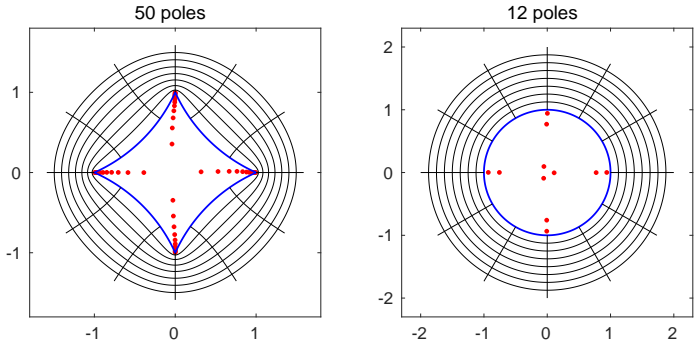}
\end{center}
\caption{\label{exterior}Exterior map for the region of
Figure~$\ref{scalloped}$ computed by a variant of Algorithm R1.\ \ The
poles of (\ref{form}) now lie interior to the boundary curve.}
\end{figure}

Just as noted at the end of the last section for smooth domains,
variants of Algorithm R1 can be developed for exterior geometries.
To map the interior of $P$ to the exterior of $S$, we change
the factor $z$ to $z^{-1}$ in (\ref{frep}) without altering the
Dirichlet problem (\ref{dirichlet}).  To map the exterior of $P$ to
the exterior of $S$, we place poles $z_j$ in (\ref{form}) interior to
$P$ rather than exterior and also change $z^j$ to $z^{-j}$.  To map
the exterior of $P$ to the interior of $S$, we make both of these
changes.  Figure~\ref{exterior} illustrates an exterior-exterior map
for the same curve $P$ as in Figure~\ref{scalloped}.  (Composing
interior and exterior maps for a given curve leads to the idea of
``conformal welding''~\cite{bishop}).)

\section{\label{smooth2}Smooth doubly connected domains}
Conformal mapping via Dirichlet problems, like the idea of
Green's functions it is based upon, does not just apply for
simply connected domains.  Following the pattern of the last few
pages, we will first introduce the mapping of smooth doubly-connected
domains in two alternative geometries, then generalize in the
next section to doubly-connected domains with corners.

Suppose $\Omega$ is the annulus bounded
by an outer Jordan curve $P_1$ and an inner Jordan curve $P_2$,
both enclosing the origin $z=0$.  Let $f$ be a conformal map of
$\Omega$ onto the circular annulus $\{w: \, \rho < |w|<1\}$
with $f(P_1) = S$ and $f(P_2) = \rho S$
for some $\rho<1$, known as the {\em conformal modulus}.
Such a map is unique up to one real degree of freedom
corresponding to a rotation.
The value of $\rho$ is unknown a priori and will
be determined as part of the solution.
Writing $f$ again in the form (\ref{frep})
shows that we must find a harmonic function $u$
in $\Omega$ satisfying
\begin{equation}
u =
\cases{
-\log |z|,  & $z\in P_1$, \cr\noalign{\vskip 5pt}
-\log |z| + \log \kern .7pt\rho,  & $z\in P_2$\kern .7pt ;
}
\label{twoeqs}
\end{equation}
see~\cite[p.~253ff\kern .7pt]{schiffer}.

Now for any value of $\rho$, not just the one we are looking
for, there will be a unique solution to the
Dirichlet problem (\ref{twoeqs})~\cite{axler}.
If $P_1$ and $P_2$ are smooth, we could approximate
this solution by functions of the form
\begin{equation}
u(z) \approx \Re \kern 1pt \Bigl( \kern 1pt\sum_{k=-n}^n c_k z^k + \alpha \log z
\Bigr)
\label{ansatz}
\end{equation}
for some real $\alpha$, with convergence as $n\to\infty$ (see the
``Logarithmic Conjugation Theorem'' of~\cite{axler}, and the final
section of the same paper for this application to doubly
connected conformal mapping).
However, if $\alpha$ is not an integer, then the analytic
function $f$ of (\ref{frep}) associated with this formula
is not single-valued.
Even if $\alpha$ is an integer, $f$ will fail to be
one-to-one in $\Omega$ if $\alpha\ne 0$.   Only with $\alpha=0$ do we
get the conformal map we seek.  And thus we find that the problem to
be solved is to find functions of the form
\begin{equation}
u(z) \approx \Re \kern -1pt \sum_{k=-n}^n c_k z^k,
\label{approx2}
\end{equation}
converging to the boundary values (\ref{twoeqs}) for\/ {\em some}
$\rho\in (0,1)$ as $n\to\infty$, and
the value of $\rho$ for which this is possible is the conformal modulus of $\Omega$.
The analogues of~(\ref{approx}) and~(\ref{vapprox}) will be
\begin{equation}
u(z) \approx \kern -3pt\sum_{k=-n}^n \kern -1pt
 \bigl( a_k \Re(z^k) - b_k \Im(z^k)\bigr),~~
v(z) \approx \kern -3pt\sum_{k=-n}^n \kern -1pt
 \bigl(a_k \Im(z^k) + b_k \Re(z^k)\bigr)~
\label{approxes}
\end{equation}
with $b_0=0$.

\begin{figure}
\vskip .2in
\begin{center}
\includegraphics[scale=.67]{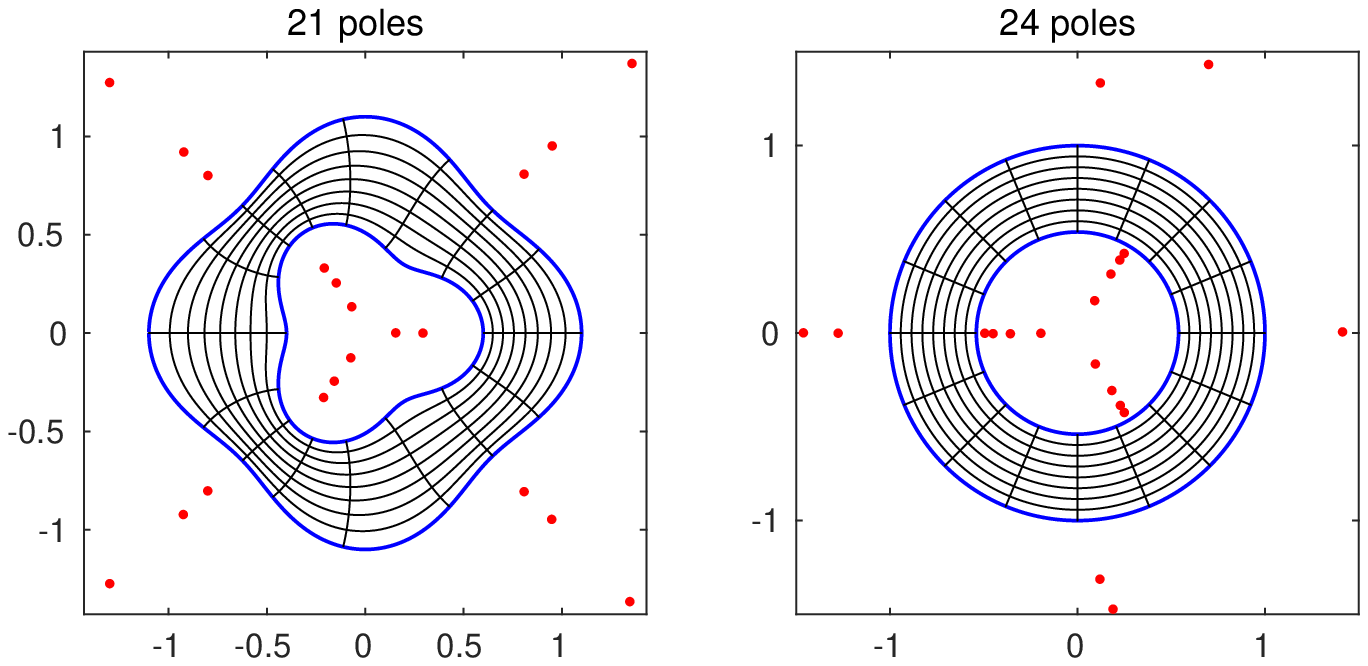}
\vskip .2in
~~\includegraphics[scale=.67]{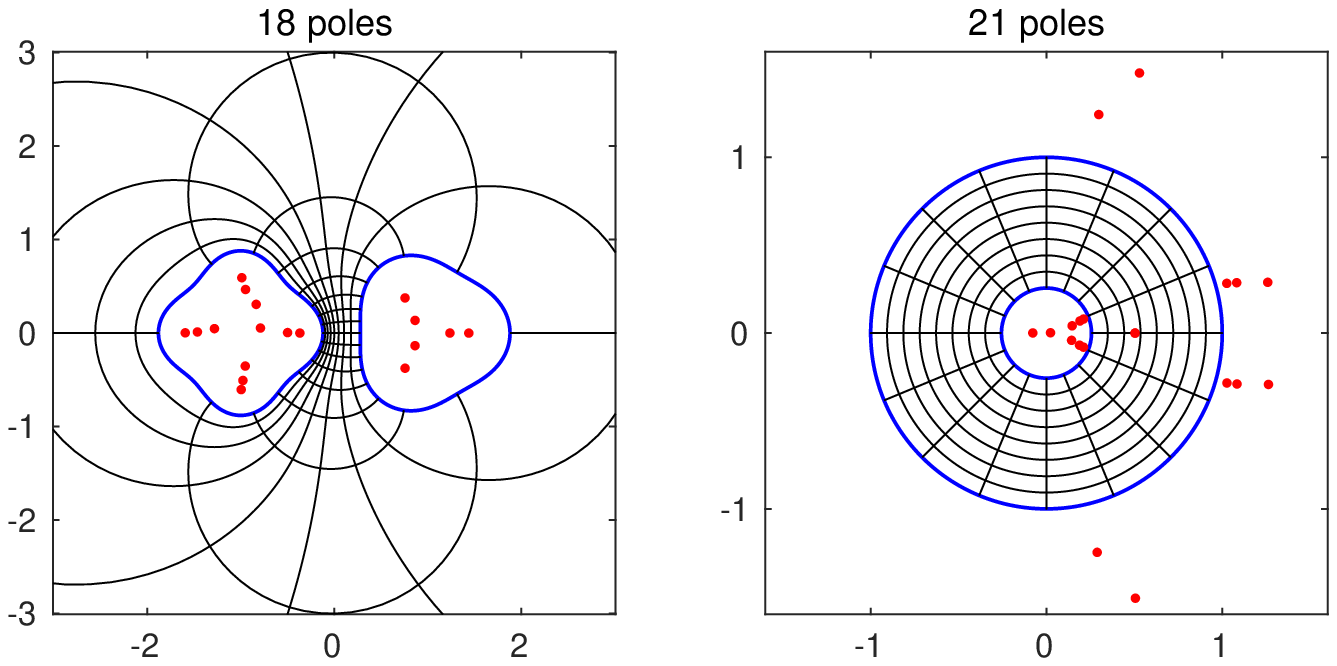}
\end{center}
\caption{\label{gaskets}Above, map of a smooth annulus computed by Algorithm P2.  Poles
of the rational representations of\/ $f$ and $f^{-1}$ appear
both inside and outside the annuli.  Below,
an unbounded smooth doubly-connected geometry mapped
via the variant\/ $(\ref{disjoint})$--$(\ref{twoeqs2})$ of Algorithm~P2.
The point $z=\infty$ is now in the interior
of $\Omega$, and its image
$f(z) = \exp(a_0) \approx 0.504$ on the right is a simple pole of $f^{-1}$.}
\end{figure}

To make (\ref{twoeqs}) and (\ref{approxes}) the basis of an
algorithm, we set up a least-squares problem whose unknowns are the
coefficients $\{a_k\}$ and $\{b_k\}$ and also the conformal modulus
$\rho$.  The number of $a_k$ and $b_k$ coefficients is $N= 4n+1$,
and the total number of unknowns is $N+1$.  The least-squares matrix
is of dimensions $2M\times (N+1)$, with the rows corresponding to the
$2M$ sample points on the boundary curves and the first $N$ columns
containing the values of $\Re(z^k)$ and $\Im (z^k)$ at these points.
The final column consists of $M$ zeros and $M$ ones, corresponding
to the appearance of $\log \kern .7pt\rho$ at one boundary curve
but not the other in (\ref{twoeqs}), and the right-hand side is the
$2M$-vector consisting of samples on the boundary of $-\log |z|$.
Here is the algorithm:

\medskip
{\em\openup 2pt

Algorithm P2.

$\bullet$ For $n = $ rounded values
of\/ $2^4, 2^{4.5}, 2^5,\dots$ until convergence or failure,

\indent\kern .19in - Set $M=8\kern .3pt n$ and
sample $-\log|z|$ in $M$ points each of\/ $P_1$ and $P_2$;

\indent\kern .19in - Solve a\/ $2M\times (N+1)$ least-squares 
problem for the coefficients\hfill\break\indent\kern .6in $a_k$ and $b_k$
of\/ $(\ref{approxes})$ and the number $\rho$;

$\bullet$ Form\/ $u$ and\/ $v$ from $(\ref{approxes})$ 
and\/ $f$ from $(\ref{frep})$\kern .7pt;

$\bullet$ Approximate $f$ by a AAA rational approximation\kern .7pt;

$\bullet$ Approximate $f^{-1}$ by another AAA rational approximation.

\par}

\medskip

\noindent 
Chebfun's {\tt conformal} command
executes this algorithm when it is called with two
boundary arcs rather than one,\footnote{{\small\sl[Note to referees:
this extension is not yet in the publicly available version of
Chebfun but will be introduced shortly.]}} and the upper half of
Figure~\ref{gaskets} shows the result computed in about
$0.2$ sec.\ by these Chebfun commands:

{
\begin{verbatim}

    z = chebfun('exp(1i*pi*z)','trig');
    C1 = z*abs(1+.1*z^4); C2 = .5*z*abs(1+.2*z^3);
    conformal(C1, C2, 'plots');      
\end{verbatim}
\par}

\noindent So far as I know, Figure~\ref{gaskets} is the first
illustration ever published of rational approximations with free
poles used to represent conformal maps of doubly-connected domains.

There is an alternative doubly-connected geometry which,
while mathematically equivalent to the one
we have just treated, may in practice be worth considering
separately.  Suppose that $\Omega$, rather than being the annular
region between two concentric curves $P_1$ and $P_2$, is the infinite region
of the extended complex plane exterior to two disjoint curves
$P_1$ and $P_2$.  If $P_1$ encloses $z=-1$ and $P_2$ encloses
$z=1$, then we can write the conformal map of $\Omega$ onto
a circular annulus in the form
\begin{equation}
f(z) = {z-1\over z+1} \kern 1pt e^{u(z) + iv(z)}
\label{disjoint}
\end{equation}
where $u$ is harmonic in $\Omega$ and satisfies, in
analogy to (\ref{twoeqs}),
\begin{equation}
u =
\cases{
-\log\Bigl|\displaystyle {z-1\over z+1}\Bigr|, 
& $z\in P_1$, \cr\noalign{\vskip 7pt}
-\log\Bigl|\displaystyle{z-1\over z+1}\Bigr| 
+ \log \kern .7pt\rho,  & $z\in P_2$
}
\label{twoeqs2}
\end{equation}
for a value of $\rho$ that again will be uniquely determined.
To fashion this into an elementary algorithm based on series, we
replace the Laurent polynomials of (\ref{approxes}) by polynomials
involving negative powers of $z-1$ and $z+1$.  The algorithm
is otherwise the same, and the second row of Figure~\ref{gaskets}
shows a computed example.  Note that $u$, $v$, and $f$ are regular at
$z=\infty$, where they take the values $a_0$, $0$, and $\exp(a_0)$,
respectively.  Thus the inverse map $f^{-1}$ has a pole in the
circular annulus at the point $\exp(a_0)$, as appears in the plot.

As with its simply connected analogue P1 of Section 2, Algorithm
P2 is far from robust, and it will quickly fail if applied to
regions much more complicated than those of Figure~\ref{gaskets}.
Superior performance can undoubtedly be achieved with integral equations,
though this has not been implemented in Chebfun.

\section{\label{corners2}Doubly connected domains with corners}
By a doubly connected domain with corners, we mean the region
$\Omega$ lying between two Jordan curves $P_1$ and $P_2$ as in the
last section, where now, each curve is a finite chain of analytic
arcs.  To conformally map such a domain, we solve the Dirichlet
problem of section~\ref{smooth2} by the lightning Laplace solver
of section~\ref{corners1}, placing poles both outside the corners
of $P_1$ and inside the corners of $P_2$.  The algorithm can be
summarized as follows.

\medskip
{\em\openup 2pt

Algorithm R2.

$\bullet$ Solve $(\ref{twoeqs})$ by
the lightning method and form $f$ from {\rm (\ref{frep})};

$\bullet$ Approximate $f$ by a AAA rational approximation\kern .7pt;

$\bullet$ Approximate $f^{-1}$ by another AAA rational approximation.

\par}

\medskip

\begin{figure}
\vskip .2in
\begin{center}
\includegraphics[scale=.67]{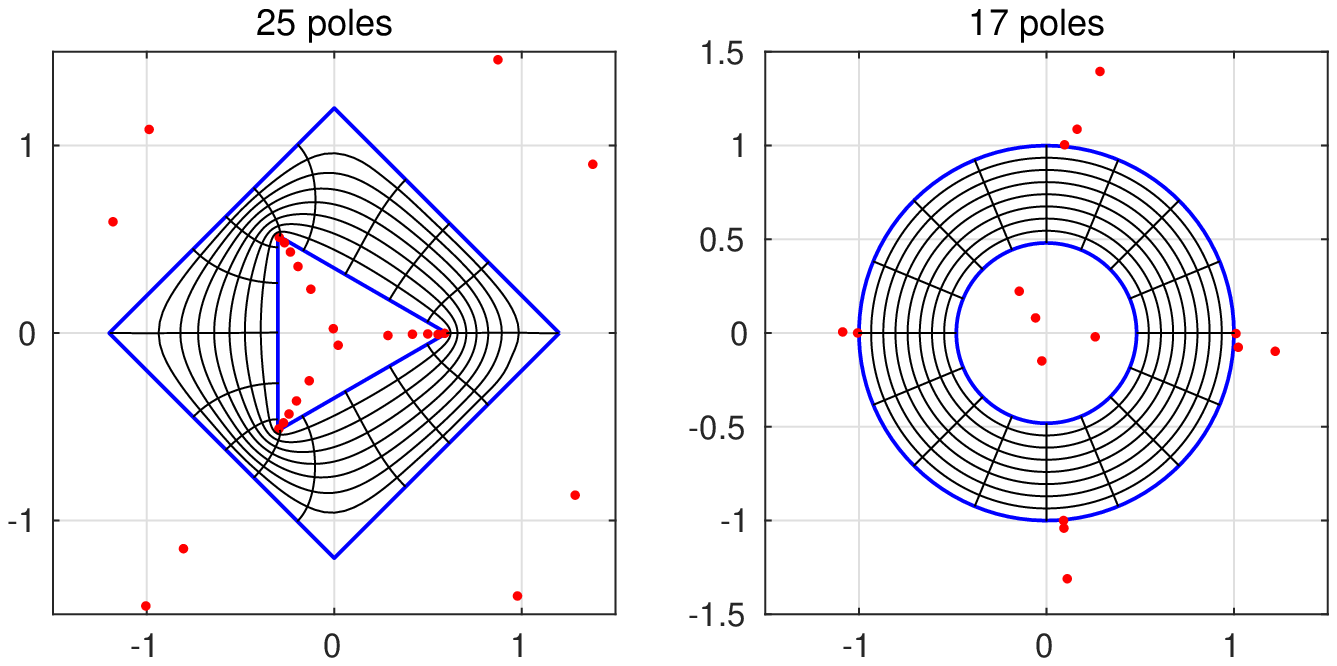}
\vskip .2in
\includegraphics[scale=.67]{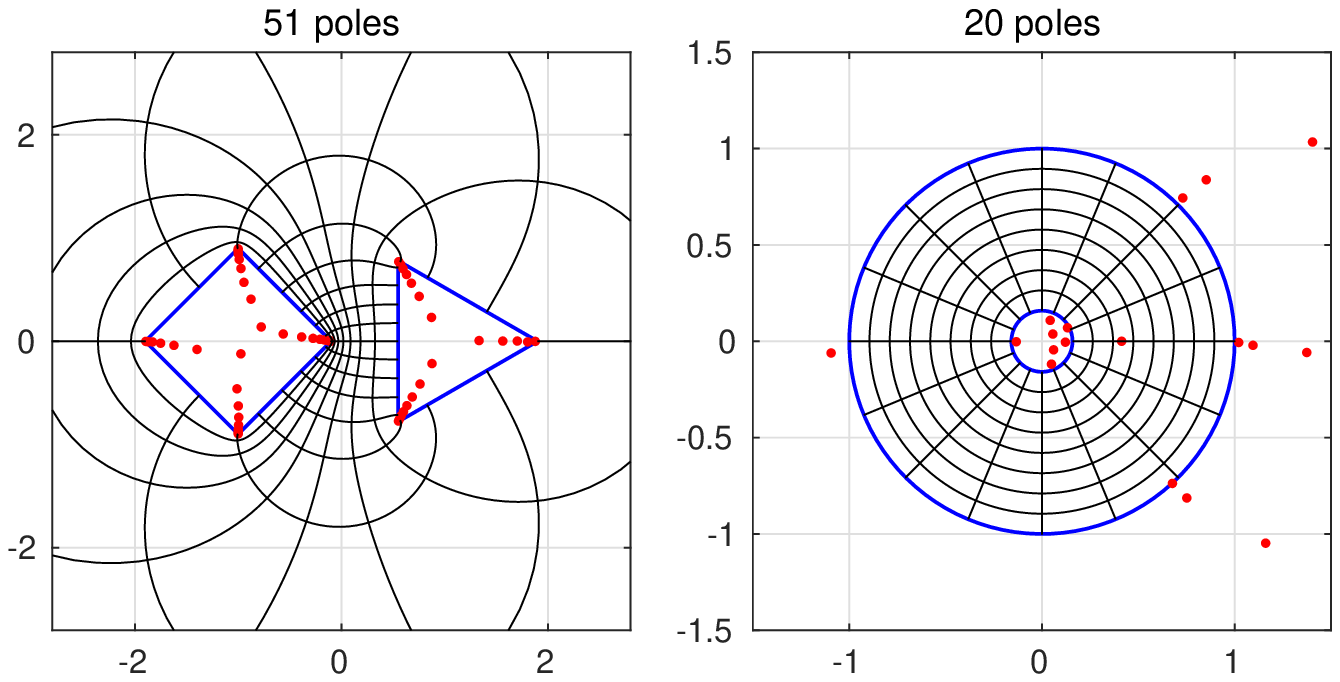}
\end{center}
\caption{\label{sharpgaskets}Doubly connected domains with corners mapped
by Algorithm R2 and its variant.  The function $f^{-1}$ in the 
lower pair again has a simple pole in the domain, shown by a dot.}
\end{figure}

\noindent The upper part of Figure~\ref{sharpgaskets} illustrates
such a map for a region that is a polygonal analogue of the region
of the first map of Figure~\ref{gaskets}.  As in the last section,
there is a straightforward analogue of Algorithm R2 for an unbounded
region exterior to two disjoint curves with corners, and this
is illustrated in the lower part of Figure~\ref{sharpgaskets}.
Such domains may be of particular interest for the calculation of
constants associated with low-rank separability, as discussed by
Beckermann, Townsend, and Wilber~\cite{bt,tw}.

For polygonal boundaries as in Figure~\ref{sharpgaskets},
there is a doubly-connected Schwarz--Christoffel
formula~\cite{dep,SCbook,henrici} which has been implemented
numerically by D\"appen~\cite{daeppen} and Hu~\cite{hu}, the latter
with a Fortran package DSCPACK.

\section{Multiply connected domains}
Domains of connectivity higher than two can also be conformally
mapped.  The Green's function idea of (\ref{twoeqs}) generalizes to
a formula with different constants $\log\rho_k$ on each boundary
component, as described in~\cite{henrici} and~\cite{schiffer}.
Various numerical methods have been developed, and some of them
are reviewed in~\cite{bds}.  The Schwarz--Christoffel formula has
generalizations that have been developed by DeLillo, Crowdy, and
others: \cite{crowdy} and \cite{dep2} are just two of many papers
in this area.  Recently Nasser has released a code PlgCirMap for
computing maps of multiply connected domains by means of a fast
implementation of an iterative method due to Koebe~\cite{nasser}.

A difficulty with 
multiply connected conformal mapping is 
the lack of an entirely natural target domain.
In the doubly connected case, the reduction to an annulus
gives a useful interpretation in almost any application.
For higher connectivity, standard choices of target
domains are regions with holes in the form of circular slits
(as delivered by the Green's function), radial slits, or disks.
Physical interpretations in such regions are more strained, however,
and it is not
always clear that much is gained by transplanting a problem such as
a partial differential equation from its native multiply-connected
geometry to one of these alternatives.

\begin{figure}
\begin{center}
\vskip .15 in
\includegraphics[scale=.67]{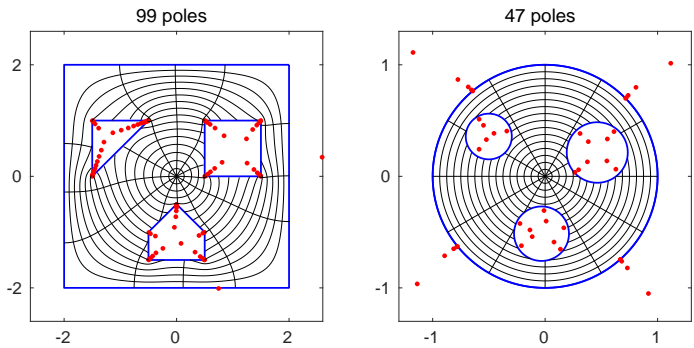}
\end{center}
\caption{\label{nasser}A multiply connected example domain
of Nasser~{\rm \cite{nasser}}.  The map is first computed
with his PlgCirMap code, and the forward and inverse maps are then
approximated by rational functions, with
poles as shown.  The rational approximations speed up
evaluations by a factor of $10$--$20$.}
\end{figure}

In view of these considerations, we have not extended our numerical
methods based on rational functions to multiply-connected regions.
If one has computed such a map by an existing method, however,
then the possibility remains of compressing it by a subsequent
rational approximation as proposed in~\cite{gt}.  Figure~\ref{nasser}
shows an example.  A 4-connected conformal map has been computed by
PlgCirMap, and AAA approximation compresses the result to a much
more efficient form.  Typical times for evaluation of the mapping
function or its inverse at each point for this example are on the
order of 200 microseconds for PlgCirMap and 10 microseconds for the
rational approximations.

\section{Numerical details}

We now mention some numerical details associated with the algorithms
presented in the previous sections.  Our implementations are for the
most part experimental, however, not claimed to be close to optimal,
and we trust that improvements in robustness and accuracy will come
in the future.

{\em Accurate data and fine sampling on the boundary.}
When failure occurs in a lightning Laplace solution, whether for
conformal mapping or more generally, the reason is usually that
the sampling on the boundary is not fine enough or the sampled
data are not accurate enough.  In particular, when poles are
exponentially clustered near a vertex, it is crucial that the
sample points be clustered on the same scale or finer.   The {\tt
laplace} software does this automatically according to a scheme
described in~\cite{lightning}, and we used the same strategy
for the doubly-connected examples of section~\ref{corners2}.
If the sample points are not dense enough, two things may go
wrong: the Laplace problem may be solved inaccurately; and, even
if it is accurate, the AAA approximation may fail, with spurious
poles appearing in regions where one expects analyticity~\cite{gt}.
(By spurious poles we mean poles of very small residue, also known as
``Froissart doublets'' since they are paired with nearby zeros that
nearly cancel their effect.)  The latter problem also arises when
the grid is fine but the data are insufficiently accurate relative
to the scale of the approximation error, so that the AAA algorithm
is effectively attempting to approximate a nonanalytic function.
Conversely, if the sample points are dense enough and the data
accurate enough, spurious poles rarely appear. 
The mechanism for this is that AAA algorithm delivers a rational
approximation with a certain kind of (linearized) optimality.
A spurious pole cannot improve an approximation much, so
it is unlikely to be present in an optimal approximation.
Success cannot be guaranteed, however, as is well known in the
theoretical literature on spurious poles in (mostly Pad\'e) rational
approximations~\cite{stahl}.

{\em Parameters of AAA approximation.}
As described in~\cite{AAA} and~\cite{lawson}, Chebfun's {\tt
aaa} code for AAA approximation has some adjustable parameters.
The overall tolerance can be specified, and we have set it
equal to the tolerance fixed for the lightning Laplace solver.
A ``cleanup'' option to remove spurious poles is invoked by default,
and we have left this in play with its tolerance set equal to the
overall AAA tolerance.  There is also an option to improve the
approximation closer to minimax by a Lawson iteration, and this
we have not invoked, because experiments suggest it is ineffective
for these problems.  All of these choices are experimental, and we
have not attempted an analysis that might justify or modify them.
As mentioned in the last paragraph, AAA approximation is far
from bulletproof.  We hope that in the next few years improvements 
to AAA will
be developed, which would have immediate impact on the robustness
of our conformal mapping algorithms.

{\em Accuracy and rate of convergence:\ theory and practice.}
The theorems of \cite{lightning} guarantee that Algorithm R1 must
converge root-exponentially, i.e.\ at a rate $C^{\kern 1pt -\sqrt n}$
for some $C>1$, so long as $\Omega$ is bounded by a finite collection
of analytic curves meeting at corners.  (We note as before that the
statements in~\cite{lightning} assume $\Omega$ is convex, but it is
believed that this condition is not actually necessary.)  Sometimes,
however, convergence stagnates at a low accuracy such as six digits.
One reason is that singularities at corners cause ill-conditioning
of the conformal map, as does the crowding phenomenon.  (A corner
of interior angle $\alpha \kern .3pt\pi$ with $\alpha<1$ will limit
the achievable number of digits to on the order of $\alpha$ times
that of machine precision.)  Another is that our method relies on
the solution of least-squares problems defined by ill-conditioned
matrices, though this does not seem to
cause much trouble in practice, a phenomenon
analyzed in the context of frames in~\cite{frames}.  

\section{Polygons and circular polygons}
In this article we have computed conformal maps by solving
Laplace problems via least-squares fitting on the boundary by
the ``lightning'' method of~\cite{lightning}, involving poles
exponentially clustered outside $\Omega$.  Thus we do not make
use of explicit analysis of singularities, as is the more familiar
approach in numerical conformal mapping.  For polygons, the famous
method based on analysis of singularities is the Schwarz--Christoffel
formula~\cite{SCbook}, which is realized numerically in the widely
used and wonderfully robust
SC Toolbox by Driscoll~\cite{toolbox}.  For circular polygons,
there is an analogous procedure based on the Schwarzian differential
equation, though it has had less impact numerically~\cite{bg,howell}.

For both polygons and circular polygons, an in-between option could
be employed that would have some advantages.  Locally near a corner
$z_j$ of a polygon with interior angle $\alpha\kern .3pt \pi$ for
some $\alpha< 2,$ the conformal map will be an analytic function of
$(z-z_j)^{1/\alpha}$ (proof by the Schwarz reflection principle).
A similar formula applies at a corner of a circular polygon,
derivable by a M\"obius transformation to a sector bounded by two
straight lines.  Thus in both of these cases, one can dispense with
exponentially clustered poles outside $\Omega$ and use a series
in terms of explicit singularities instead.  At the corner of
the polygon, for example, the terms will be $(z-z_j)^{1/\alpha},
(z-z_j)^{2/\alpha}, (z-z_j)^{3/\alpha},\dots .$ By constructing
the matrix $A$ from columns corresponding to such terms, one can
set up a least-squares problem just as before, but making use
of a very different set of basis functions, and the convergence
will be exponential rather than root-exponential.  This would be a
numerical method midway between the traditional Schwarz--Christoffel
or Schwarzian differential approaches and the fully general one
we have focused on in this paper.  (A disadvantage is that the
inverse map would have to be treated separately, if AAA
rational approximation is to be avoided.)  We have not pursued this idea,
but proofs of concept for L-shaped and circular L-shaped regions,
with short Matlab codes, can be found in~\cite{example}.

One might suppose that these observations indicate that there is no
need for lightning Laplace solvers in general when the domain is a polygon or
circular polygon.  However, this is not so, for the Laplace problems
that arise in conformal mapping are special ones, with corner
singularities determined only by the local geometry~\cite{lehman}.
More general Laplace problems often have different boundary data
on different sides of a domain, and in such a case, the corner
singularities may be complicated even when the geometry is
simple~\cite{lehman2}.

\section{Discussion}

The algorithms proposed in this paper combine two phases.
First, a Laplace problem is solved by least-squares collocation
on the boundary.  This is essentially the idea of the Method of Fundamental
Solutions, which has been pursued by Amano and his coauthors for
conformal mapping~\cite{amano,amano2} under the
name of the charge simulation method.  However, our approach differs in the
use of poles exponentially clustered near corner singularities to
get root-exponential convergence, as proposed in~\cite{lightning}
(or in Section 8, the use of explicit singular terms), and
in the use of rational functions (dipoles) rather than
point charges (monopoles).  Second,
the resulting map and its inverse are both approximated by rational
functions, as proposed in~\cite{gt}.

The idea of solving Laplace problems by fitting Dirichlet data on
the boundary is an old
one, and early discussions of note are by
Walsh~\cite{walsh}, Curtiss~\cite{curtiss},
Collatz~\cite[sec.~V.4.2]{collatz},
and Gaier~\cite[sec.~4.1]{gaier},
the last of these specifically for conformal mapping.
However, these authors formulated systems of equations
rather than least-squares problems, making the method fragile and
impractical.\footnote{Gaier:
``Since the configuration of the maximal system of points on
the boundary is unknown (except for circles and ellipses), this
method is mainly of theoretical interest.''
Curtis: ``The success of this
method is highly sensitive to the correct placement of the
points.'' Collatz: ``The choice of
collocation points is a matter of some uncertainty.''}
As soon as one moves from
square matrices and systems of equations to rectangular matrices
and least-squares, the distribution
of sample points largely ceases to be a problem.

Our methods are based on a minimum of geometrical analysis,
and we do not claim they are optimal for any particular
problem.  For mapping polygons, in particular, the SC Toolbox is tried
and true and able
to handle quite extreme geometries without difficulty~\cite{toolbox}.
Nevertheless, the approach to conformal mapping presented here
is very general and usually very fast, enabling
typical problems to be solved in on the order of a second of computer
time, with subsequent evaluations of both the map and its inverse in
microseconds.  Experimental codes can be found at~\cite{lightcode}.
Extensions to problems such as domains with slits, unbounded domains, or domains
bounded by arcs other than straight lines or arcs of circles should
be straightforward, though they have not yet been explored.

\begin{acknowledgements}I have benefited from helpful comments
of Toby Driscoll, Abi Gopal, and Yuji Nakatsukasa.
\end{acknowledgements}

\end{document}